\documentclass[12 pt]{amsart}

\usepackage{amscd, amssymb}

\newtheorem{pr}{Proposition}
\newtheorem{lm}{Lemma}
\newtheorem{tm}{Theorem}
\newtheorem{cor}{Corollary}
\newcommand{\bpf}{\noindent {\em Proof.} }
\newcommand{\epf}{\qed \vspace{+10pt}}

\newcommand{\proj}{\mathbf P}

\newcommand{\rarr}{\rightarrow}
\newcommand{\oh}{{\mathcal{O}}}
\newcommand{\com}{\mathbb{C}}

\newcommand{\Z}{\mathbb{Z}}

\newcommand{\br}{{br}}
\newcommand{\coh}{_{coh}}

\newcommand{\sym}{{\text{Sym}}}
\newcommand{\comm}{^{\bullet}}
\newcommand{\KK}{{\mathcal{K}}}

\begin{document}
\title{Stable maps and branch divisors}
\author{B. Fantechi and R. Pandharipande}
\date{18 May  1999}
\maketitle

\pagestyle{plain}
\setcounter{section}{-1}
\section{\bf{Introduction}}

Let $f:X\to Y$ be a surjection of nonsingular projective varieties of the 
same dimension. The ramification divisor $R$ of $f$ on $X$ is defined 
by requiring the sequence
\begin{equation}
\label{basiccase}
0 \rarr {f^*\omega_Y}\rarr {\omega_X} \rarr {\omega_X|_R} \rarr 0
\end{equation}
to be exact. The branch divisor $\br(f)$ on $Y$ is then defined by
pushing forward: $\br(f)=f_*(R)$.
The support of $\br(f)$ is the locus of 
points $y\in Y$ 
such that $f$ is not \'etale in any neighborhood of $f^{-1}(y)$.

If $f:C \rarr D$ is a degree $d$ map of
nonsingular curves,
then $br(f)$ is a divisor on $D$ of degree 
$$r=2g(C)-2-d(2g(D)-2)$$
by the Riemann-Hurwitz formula.
Let $M_{g}(D,d)$ the moduli stack of degree $d$ maps from 
nonsingular genus $g=g(C)$ curves to $D$.
The branch divisor yields a morphism
of Deligne-Mumford stacks
\begin{equation}
\label{oooo}
\gamma: M_{g}(D,d) \rarr \text{Sym}^{r}(D).
\end{equation}
For moduli points $[f:C \rarr D] \in M_{g}(D,d)$, 
$\gamma([f])= br(f)$.
A natural extension of $\gamma$ to the compactification
by stable maps $$M_g(D,d) \subset  
\overline{M}_{g}(D,d)$$ is the main result of the
paper.

Consider first the following situation. Let
$f:X \rarr Y$ be a projective morphism of
$S$-schemes where:
\begin{enumerate}
\item[(i)] $X$ is a
local complete intersection over $S$ of relative dimension $n$.
\item[(ii)] $Y$ is smooth over $S$ of relative dimension $n$.
\item[(iii)] All geometric fibers of $X$ over $S$
are reduced.
\end{enumerate}
Under these conditions, 
a functorial
 relative Cartier divisor $\br(f)$ on $Y$ over $S$ is 
constructed in Section \ref{branchcon}.
The divisor $\br(f)$ is 
supported on the locus of 
points $y\in Y$ such 
that $f$ is not \'etale in any neighborhood of $f^{-1}(y)$.
In this generality, $\br(f)$ need not be an effective Cartier divisor. 
However, $\br(f)$ is invariant 
under  base change 
and coincides with the branch divisor  defined by (\ref{basiccase})
when  $X\to S$ is smooth and every component of $X$
dominates one of $Y$.

The branch divisor $\br(f)$ is constructed by studying the complex
\begin{equation}
\label{cmplx}
Rf_*[f^*\omega_{Y/S}\to \omega_{X/S}],
\end{equation}
well-defined up to isomorphism in $D^-\coh(Y)$. By
generalizing to complexes a classical construction 
of Mumford for sheaves ([Mu], \S 5.3), 
we can associate to (\ref{cmplx}) a Cartier divisor on $Y$.
Section \ref{pertor} contains the required generalization of
Mumford's results.

In Section \ref{exten}, we apply our branch divisor construction
to the universal family:
$$F: {\mathcal{C}} \rarr D \times \overline{M}_{g}(D,d)$$
over the moduli stack of
stable maps $\overline{M}_{g}(D,d)$ for $d>0$.
Certainly this universal family (as a Deligne-Mumford stack)
satisfies conditions 
(i-iii). 
It is shown $\br(F)$ in this case is an {\em effective}
relative Cartier divisor on $D \times \overline{M}_{g}(D,d)$
of relative degree $r$.
The branch divisor $\br(F)$ then yields a canonical
morphism  
\begin{equation}
\label{exxx}
\gamma: \overline{M}_{g}(D,d)\rarr  \sym^{r}(D)
\end{equation} 
extending (\ref{oooo}).

The morphism $\gamma$ has an appealing point theoretic
description on the boundary of $\overline{M}_{g}(D,d)$.
Let $[f: C \rarr D]$ be a moduli point where
$C$ is a singular curve. Let $N\subset C$ be the
cycle of nodes of $C$.
Let $\nu: \tilde{C} \rarr C$
be the normalization of $C$.
Let $A_1, \ldots, A_a$ be the components of $\tilde{C}$
which dominate $D$, and let $\{a_i: A_i \rarr D\}$
denote the natural maps. 
As $a_i$ is a surjective map between nonsingular
curves, the branch divisor $\br(a_i)$ is defined by (\ref{basiccase}). 
Let $B_1, \ldots, B_b$ be
the components of $\tilde{C}$ contracted over $D$, and
let $f(B_j)=p_j\in D$.
We prove the formula:
\begin{equation}
\label{ptwise}
\gamma([f])=\br(f)= \sum_{i} \br(a_i) +
\sum_j (2g(B_j)-2)[p_j] + 2f_{*}(N).
\end{equation}
It is easy to see that formula (\ref{ptwise}) associates an
effective divisor of degree $r$ on $D$ to every
moduli point $[f]$. However, the construction of 
$\gamma$ as a scheme-theoretic morphism requires
the relative branch divisor results over arbitrary
reducible, nonreduced bases $S$.

In Section \ref{hurnum}, the morphism
$\gamma$ is used to study the
classical Hurwitz numbers $H_{g,d}$
via Gromov-Witten theory.
$H_{g,d}$ is 
the number of nonsingular, genus $g$ 
curves expressible as  $d$-sheeted covers
of $\proj^1$ with a {\em fixed} general branch divisor.
The Hurwitz numbers were first computed in [Hu] by 
combinatorical techniques. A simple analysis of the
moduli space of stable maps to $\proj^1$ shows:
\begin{equation}
\label{hurw}
H_{g,d}= \int_{ [\overline{M}_{g}(\proj^1,d)]^{vir}} 
\gamma^*(\xi^{2g-2+2d}),
\end{equation}
where $\xi$ is the hyperplane class on 
$\sym^{2g-2+2d}(\proj^1)=\proj^{2g-2+2d}$.
It is then possible to directly evaluate the integral (\ref{hurw})
using the virtual localization formula [GrP] to obtain
a Hodge integral expression for the Hurwitz numbers:
\begin{equation}
\label{hodgehur}
H_{g,d}= \frac{(2g-2+2d)!}{d!} \int_{\overline{M}_{g,d}}
\frac{1-\lambda_1+\lambda_2 -\lambda_3 + \ldots + (-1)^g \lambda_g}
{\prod_{i=1}^d (1-\psi_i)},
\end{equation}
for $(g,d) \neq (0,1), (0,2)$.
The integral on the right is taken over the moduli
space of pointed stable curves $\overline{M}_{g,d}$. 
The classes $\psi_i$ and $\lambda_j$
are the cotangent line classes and the Chern classes of the
Hodge bundle respectively. The values  $H_{0,1}=1$ and $H_{0,2}=1/2$
are degenerate  cases from the point of view of the right side
of (\ref{hodgehur}).

A proof of formula (\ref{hodgehur}) has been announced
independently by Ekedahl, Lando, Shapiro, and
Vainshtein using very different methods [ELSV].
In fact,
the formula of [ELSV] also accounts for particular non-simply
branched cases which appear again to be equal to
vertex integrals in the
virtual localization formula. However, the connection in the
non-simply branched cases is not clear. 

We thank
T. Graber and R. Vakil for many discussions about
Hurwitz numbers. 
The authors thank the Mittag-Leffler Institute where
this research was partially carried out. 
The second author was
partially supported by National Science
Foundation grant DMS-9801574.

\section{\bf{Perfect torsion complexes}}
\label{pertor}
\subsection{Cartier divisors}
The base field $\com$ of complex numbers will
be fixed for the entire paper. However, all the results
of Sections 1.1 - 2.2 are valid over any algebraically
closed base field. The characteristic 0 condition
is required for generic smoothness in the
construction of the branch divisor.

Let $A$ be an algebra of finite type over 
$\com$. 
Let $S\subset A$ be the multiplicative system of 
elements which are not zero divisors. 
Recall, the set of zero divisors of $A$ equals the union
of all associated primes of $A$ ([Ma], p.50). 
A prime ideal $\mathfrak{p}\subset A$ is
{\em depth 0} if all non-units of 
$A_{\mathfrak{p}}$ are zero divisors.
The associated primes of $A$ are exactly the depth
0 primes ([Ma], p.102).
Let $K(A)= S^{-1}(A)$ be the total quotient ring of
$A$. It is easy to check for $f\in A$,
$K(A_f)= K(A)_f$.

Let $X$ be a scheme (always taken here to be
quasi-projective over $\com$). 
We distinguish the points of $X$ (integral subschemes)
from the geometric points of $X$ ($Spec(\com)$ subschemes).
Let $\KK$ be the sheaf of rings on $X$ defined
by associating $K(A_i)$ to the  basis of all affine
open sets $Spec(A_i)$ of the Zariski topology of $X$.
The equality 
$$\Gamma(Spec(A_i), \KK)= K(A_i)$$ follows from
 the property $K(A_f)=K(A)_f$.
Let $\KK^*$ denote the sheaf of invertible elements of $\KK$.
A Cartier divisor is an element of $\Gamma(X, \KK^*/\oh^*)$.
This discussion follows Hartshorne ([Ha1], \S II.6).  

A Cartier divisor is defined by the data
$\{(f_i, W_i)\}$ where the open sets 
$W_i=Spec(R_i)$ cover $X$ and
$$f_i \in K(R_i)^*, \ \ \ f_i/f_j \in \Gamma(W_i \cap
W_j, \oh^*).$$
A Cartier divisor $D$ is {\em effective} if
there exist defining data as above satisfying $f_i \in R_i
\subset  K(R_i)$. An effective Cartier divisor
naturally defines a locally free ideal sheaf of $\oh_X$.
\begin{lm}
Let $U\subset X$ be an open set containing all
depth 0 points of $X$.
Let $f\in \Gamma(U, \oh_U^*)$. Then, $f$
defines a canonical element of $\Gamma(X, \KK^*)$.
\end{lm}

\bpf
Let $Z=U^c \subset X$.
Let $\{ W_i =Spec(R_i)\}$ be an open affine cover of $X$.
Let $U_i=U\cap W_i$, $Z_i=Z\cap W_i$, and $f_i=f|_{U_i}$. 
Let $I\subset R_i$ be the radical ideal determined
by closed set
$Z_i$.
Since $Z_i$ contains no depth 0 points,
$I$ must contain a element $x$ of $R_i$ which is not
a zero divisor. Since $Spec((R_i)_{x}) \subset U_i$, we
see $f_i$ is naturally an element of
$(R_i)^*_{x}$. As $K(R_i)$ is obtained from
$(R_i)_x$ by further localization, $f_i$ yields a
canonical element of $K(R_i)^*$.
These local sections over $W_i$ patch to yield
a canonical element of $\Gamma(X, \KK^*)$.
\epf

\subsection{The divisor construction (local)}
We recall here a construction of Mumford ([Mu], \S 5.3). 
For our general branch
divisor construction, we must extend these results from
sheaves to complexes.

Let $D^-\coh(X)$ be the derived category of bounded (from above)
complexes 
of quasi-coherent 
$\oh_X$-modules with coherent cohomology on a scheme $X$.
We will identify a sheaf with a complex in degree zero;
we will identify a morphism with a complex in degrees $[-1,0]$.
By convention, {\em free} and 
{\em locally free} sheaves will have finite rank.
An object $E\comm$ of $D^-\coh(X)$ is {\em perfect} 
if it is locally isomorphic to a finite complex of  
locally free sheaves.
$E\comm$ is  {\em torsion} if for 
all $i\in\Z$
the support of $H^i(E\comm)$ does not contain 
any point of depth zero of $X$.

Let $E\comm=[E^a\to E^{a+1}\to\ldots\to E^b]$ be a finite 
complex of free sheaves on $X$, and let $rank(E^i)=r_i$. 
Let $$\Lambda(E\comm)=\bigotimes\limits_{i=a}^b
(\Lambda^{r_i}E^i)^{(-1)^i}.$$
Following [Mu], a {\em choice} of an explicit isomorphism 
$E^i=\oh_X^{r_i}$ for each $i$ yields 
an isomorphism $$\psi :\Lambda(E\comm)\to \oh_X.$$ Because of 
the choice of the trivializations of $E^i$,
$\psi$ is determined only up to multiplication 
by a section of $\oh_X^*$.
However, if $E\comm$ is exact, there is a canonical isomorphism 
$$\kappa:\Lambda(E\comm)\to \oh_X.$$
These isomorphisms $\psi$ and $\kappa$ will together determine
a Cartier divisor in the torsion case.

Let $E\comm=[E^a\to E^{a+1}\to\ldots\to E^b]$ be a 
finite torsion complex of free sheaves on $X$. 
We define the associated Cartier divisor $div(E\comm)$ following
 [Mu]. 
Let $U$ be the complement of the union of the supports 
of $H^i(E\comm)$.
On $U$ there is a canonical 
isomorphism $$\kappa_U:\Lambda(E\comm)|_U\rarr\oh_U.$$ 
Let $\psi :\Lambda(E\comm)\to \oh_X$ be an isomorphism 
defined by trivializations
as above.
Then, $\psi_U\circ(\kappa|_U)^{-1}$ is a section 
of $\oh^*_U$. As $U$ contains all points of depth zero of $X$, 
we obtain a unique section $f$ of $\KK^*_X$ by Lemma 1.
Different 
trivializations of $E^i$ over $X$
change the isomorphism $\psi$ 
by multiplication by an element of $\oh_X^*$.
 Hence, $f$ yields a well-defined section of 
$\KK^*_X/\oh_X^*$. 
Let $div(E\comm)$ denote this
canonically associated Cartier divisor on $X$.
Note if $E\comm$ is an exact finite complex of free sheaves,
then $div(E\comm)$ is zero.

\subsection{The divisor construction (global)}

Let $\phi:E\comm\to F\comm$ be a chain map of complexes.
The mapping cone of $\phi$ is the complex $M(\phi)\comm$ with
sheaves
$M(\phi)^i=E^{i+1}\oplus F^i$ and differentials
$M(\phi)^{i-1}\to M(\phi)^{i}$ determined
by $(e,f)\mapsto (de,df+(-1)^i\phi(e))$ (where $d$ denotes
differentials on $E\comm$ and $F\comm$).
Note  there are natural morphisms of complexes 
\begin{equation}
\label{llll}
E\comm\to F\comm\to M(\phi)\comm \to E\comm[1]
\end{equation}
where $F^i\to M(\phi)^i$ is given by 
$f\mapsto (0,f)$ and 
$M(\phi)^i\to E^{i+1}$ is defined by $(e,f)\mapsto e$.
Any sequence of morphisms $E\to F\to G\to E[1]$ in $D^-\coh(X)$ 
which is isomorphic to
(\ref{llll}) in $D^-\coh(X)$ is called a {\em distinguished triangle}.

The morphisms (\ref{llll}) induce 
a long exact sequence of cohomology $$
\dots \to H^i(E)\to H^i(F)\to H^i(M(\phi))\to H^{i+1}(E)\to \dots$$
In particular, if $\phi$ is a quasi-isomorphism, then $M(\phi)
\comm$ is exact.

\begin{lm} 
\label{connn}
Let $E\comm$ and $F\comm$ be  finite torsion complexes
of free sheaves, and let $\phi:E\comm \to F\comm$ be a chain map. 
Then, the mapping cone $G\comm$ of $\phi$ is also a 
finite torsion complex of free sheaves, 
and $$div(F\comm)=div(E\comm)+div(G\comm).$$
\end{lm}
\bpf
$G\comm$ is certainly 
 a finite complex of free sheaves.
Let $Z\subset X$ be 
the union of the supports of the cohomology
sheaves of $E\comm$ and $F\comm$. 
As the latter supports do not contains points of depth zero, 
neither does $Z$. Let $U= Z^c$.  Both $E\comm$ and $F\comm$ are
exact on $U$, so $\phi|_U$ is a quasi-isomorphism and 
$G\comm|_U$ is also exact. Hence, $G\comm$ is torsion.

There is a canonical isomorphism of 
$\Lambda(F\comm)$ with $\Lambda(E\comm)\otimes\Lambda(G\comm)$, 
which proves the lemma. 
\epf

\begin{cor} 
\label{ccc}
Let $E_1\comm$ and $E_2\comm$ be finite 
torsion complexes of free sheaves. If they are isomorphic in $D^-\coh(X)$, 
then the induced Cartier
divisors $div(E_1\comm)$ and $div(E_2\comm)$ are equal. 
\end{cor}
\bpf
If $E_1\comm$ and $E_2\comm$ are isomorphic in $D^-\coh(X)$, 
then there exists
an object $L\comm\in D^-\coh(X)$ and 
chain maps $L\comm\to E_i\comm$ which are
quasi-isomorphisms.
We may prove  
the Corollary locally on $X$. Locally, 
we can find a free 
complex $F\comm$ with a chain map 
$F\comm\to L\comm$ which is a quasi-isomorphism ([Ha1], Lemma 12.3).
As $E_i\comm$ are finite and free, $F\comm$ may be cut-off from
below to yield a finite and free complex with
quasi-isomorphisms: $F\comm_{cut} \rarr E_i\comm$.

It is therefore enough to prove the Corollary 
in case there exists
a quasi-isomorphism $\phi:E_1\comm\to E_2\comm$, 
but then it follows from
Lemma \ref{connn}.
\epf

Let $E\comm$ be a perfect torsion complex on $X$.
As $E\comm$ is locally a finite torsion complex of free
sheaves, Cartier divisors may be associated locally 
to $E\comm$ via local trivializations and the 
construction of Section 1.2. 
By Corollary \ref{ccc}, these locally associated divisors agree
and define a canonical
Cartier divisor $div(E\comm)$ on $X$.

\begin{pr}
\label{mainprop}
Let $E\comm$ be a perfect torsion complex on $X$.
Then $div(E\comm)$ satisfies the following properties:
\begin{enumerate}
\item[(i)]
$div(E\comm)$ depends only on  
the isomorphism class of $E\comm$ in $D^-\coh(X)$.
\item[(ii)]
 If $F$ is a coherent torsion sheaf on $X$ admitting locally a finite free
resolution, then $div(F)$  is the divisor  constructed 
in [Mu]. Moreover, $div(F)$ is an {\em effective} Cartier
divisor.
\item[(iii)]
If $D$ is an effective Cartier divisor
in $X$, $div(\oh_D)=D$.
\item[(iv)]
The divisor is additive for distinguished triangles.
\item[(v)]
If $f:X'\to X$ is a base change, such that $f^*E\comm$ is torsion, then
$f^*(div(E\comm))$ is a Cartier divisor. Moreover, in this case 
$$div(f^*(E\comm))=f^*(div(E\comm)).$$
\item[(vi)]
$div(E\comm[-1])=-div(E\comm).$
\item[(vii)]
If $L$ is a line bundle on $X$, $div(E\comm)=div(E\comm\otimes L)$.
\end{enumerate}
\end{pr}

\bpf
For the most part, these properties are simple
consequences of the construction.
Property
(i) follows immediately from local considerations and  Corollary
\ref{ccc}.
The equivalence with Mumford's construction (ii) is true by definition.
The effectivity of $div(F)$ is a subtle issue proven in [Mu].
An easy computation 
using the isomorphism between
 $[\oh_D]$ and $[\oh_X(-D)\rarr \oh_X]$ in
$D^-\coh(X)$ proves (iii). Lemma (\ref{connn}) and
local analysis together imply (iv).
Property (v) may be checked locally on $X$ and $Y$ where
the divisor construction is seen to be compatible with
the definition of the pull-back of Cartier divisors.
Properties (vi) and (vii)
are trivial consequences of the definitions.
Property (vi) shows $div(E\comm)$ is not an effective
Cartier divisor for all perfect torsion complexes.
\epf

The following  example will be required.
Let $X$ be a projective scheme, and let $Y$ be
a nonsingular curve.
Let $f:X\rarr Y$ be a constant morphism with image $y\in Y$. 
\begin{lm}
\label{smashpoint}
For any coherent sheaf $F$ on $X$,
$Rf_*(F)$ is a perfect torsion complex in $D^-\coh(Y)$, and
$div(Rf_*(F))=\chi(F)[y]$.
\end{lm}
\bpf
$Rf_*(F)$ defines a complex in $D^-\coh(Y)$ with
coherent cohomology, nonzero in finitely many degrees. 
By the nonsingularity of $Y$,
$Rf_*(F)$ is perfect. That $Rf_*(F)$ is torsion is
clear. As $Rf_*(F)$ is exact outside of $y$, $div(Rf_*(F))$
is a multiple of the point $[y]$. The Lemma then follows
from a local calculation.
\epf

\subsection{Torsion criterion}
Let $q:Y \rarr S$ be a smooth morphism with irreducible
fibers. Let $Ass(Y)$ and $Ass(S)$ be the sets of
depth 0 points
of the schemes $Y$ and $S$ respectively. 
A point ${\mathfrak{p}}$ of $S$ corresponds to an integral subscheme
$V_{\mathfrak{p}} \subset S$. Since $q$ is smooth with irreducible fibers,
$q^{-1}(V_{\mathfrak{p}})$ 
is an integral subscheme of $Y$ determining
a point ${\mathfrak{q}}$ of $Y$. Let $\iota({\mathfrak{p}})
={\mathfrak{q}}$.
\begin{lm} $\iota( Ass(S)) = Ass(Y)$.
\label{tqtqtq}
\end{lm}
\bpf
The Lemma may be checked locally on $Y$ and $S$, so we may take
$Y=Spec(B)$ and $S=Spec(A)$. Since $q$ is smooth, $q$ is
flat. 
If $M$ is a Noetherian $R$-module, let $Ass_R(M)$ denote
the set of primes of $R$ associated to $M$.
An algebraic result from Bourbaki is now required
(also [Ma], Theorem 12):
\begin{equation}
\label{dsdsds}
Ass_B(B)= \bigcup_{{\mathfrak{p}} \in Ass_A(A)} 
Ass_B(B/{\mathfrak{p}}B).
\end{equation}
As discussed above, ${\mathfrak{p}B} \subset B$ is a
prime ideal. Hence $Ass_B(B/{\mathfrak{p}}B)= 
\{ {\mathfrak{p}}B \}$.
Moreover, $\iota({\mathfrak{p}})= {\mathfrak{p}}B$ by
definition.
\epf

Let $E\comm$ be a
perfect object of $D^-\coh(Y)$. We will require the following
criterion for torsion.
\begin{lm} 
\label{qtqtqt}
Let $q:Y \rarr S$ be a smooth morphism with irreducible fibers. 
If for every geometric point
$s\in S$, the
complex $i_s^*(E\comm)$ is torsion on $Y_s$
(where $i_s:Y_s\rarr  Y$ is the inclusion), 
then $E\comm$ is torsion on $Y$.
\end{lm}

\bpf 
We again may take $Y=Spec(B)$ and $S=Spec(A)$.
Let ${\mathfrak{q}}= \iota({\mathfrak{p}})$ be
a depth 0 point of $Y$. 
By Lemma \ref{tqtqtq}, all depth 0 points of $Y$ may be so expressed.
Let $y\in V_{{\mathfrak{q}}}$
be a geometric point of $Y$ with $s=q(y)$ satisfying:
$i_s^*(E\comm)$ has cohomology supported away from
$y$ in $Y_s$. Such a $y$ can be found since
$V_{{\mathfrak{q}}}$ contains fibers of $q$.
As $E\comm$ is perfect on $Y$, we can take a finite locally
free representative
$$E\comm= [E^a \rarr E^{a+1} \rarr \cdots \rarr E^b]$$
locally at $y\in Y$. Since the fiber sequence 
$$0 \rarr E^a_y \rarr E^{a+1}_y \rarr \cdots \rarr E^b_y\rarr 0$$
is exact by the torsion condition on $i_s^*(E\comm)$,
$E\comm$ is exact in a Zariski neighborhood of $y$ in $Y$.
In particular, the point $y$ does not
lie in the support of the cohomology 
of $E\comm$ on $Y$. Since
$y$ is in the closure of the point ${\mathfrak{q}}$,
we see ${\mathfrak{q}}$ does not lie in the
cohomology support. 
\epf

We first note $i_s^*(E\comm)$ is the pull-back in
the derived category. For a complex of free objects
(or, more generally a complex of $S$-flat objects), this pull-back
is determined by the simple pull-back of sheaves. Second,
we note the irreducibility hypothesis  on the fibers
of $q:Y \rarr S$ can be easily removed in Lemma \ref{qtqtqt}
by generalizing Lemma \ref{tqtqtq} slightly. We leave the details
to the reader.

\section{\bf{Branch divisors}}
\label{branchcon}
\subsection{Notation}
\label{sss}
Let $X$, $Y$, and $S$ be schemes. Let
$$p: X \rarr S, \ \ q: Y\rarr S$$
be morphisms satisfying:
\begin{enumerate}
\item[(i)] $X$ is a
local complete intersection over $S$ of relative dimension $n$.
\item[(ii)] $Y$ is smooth over $S$ of relative dimension $n$.
\item[(iii)] All geometric fibers of $X$ over $S$ are reduced.
\end{enumerate} 
Let $f: X \rarr Y$ be a projective morphism over $S$.
This data will be fixed for the entire section. 
We will construct a relative Cartier divisor $br(f)$ on $Y$ generalizing
the standard branch divisor. 

\subsection{Direct images}

We review here the natural map $$Rf_*: D^-\coh(X) \rarr D^-\coh(Y)$$
obtained from direct images.
Let ${\mathcal{U}}$ be an $f$-relative Cech cover of $X$
(over every affine open in $Y$, ${\mathcal{U}}$ restricts to
a usual Cech covering). For any quasi-coherent sheaf
$E$ on $X$, let  $C\comm({\mathcal{U}},E)$ be the associated
Cech complex of quasi-coherent sheaves on $Y$.
Let $E\comm$ be an object of $D^-\coh(X)$.
Then, $Rf_*(E\comm)$ is defined
to be the simple complex on $Y$
obtained from the double complex
$C^p({\mathcal{U}}, E^q)$. The complex $Rf_*(E\comm)$ is
certainly bounded from above. Moreover, the cohomology of 
$Rf_*(E\comm)$ may be computed by a spectral sequence with $E_2$ term
$R^pf_*(H^q(E\comm))$.
Since, $R^pf_*(H^q(E\comm))$ is a
grid of coherent sheaves on $Y$ with only finitely many
objects on each line of slope $-1$, the cohomology of $Rf_*(E\comm)$
is coherent.
Hence, $Rf_*(E\comm)$ defines an element of $D^-\coh(Y)$.
To show this construction is well-defined in the
derived category, see [Ha2].
\begin{lm}
$Rf_*:D^-\coh(X)\to D^-\coh(Y)$
carries perfect complexes to perfect complexes.
\label{perfperf}
\end{lm}

\bpf
The statement is local, so we assume $Y$ is affine.
Since $f$ is projective and $E\comm$ is perfect, we can assume
$E\comm$ is a finite complex of locally free sheaves globally on $X$. 
By Lemma 5.8
of [Mu], each of the Cech sheaves $C^p({\mathcal{U}},E^q)$
has finite $Tor$-dimension and 
hence admits a finite flat resolution by quasi-coherent sheaves
on $Y$.
Therefore $Rf_*(E\comm)$ is isomorphic in the derived
category to a finite complex of quasi-coherent 
flat sheaves and hence is $Tor$-finite.
As $Rf_*(E\comm)$ is bounded from above and has coherent 
cohomology, we can construct an isomorphic complex of 
locally free sheaves, indexed in $(-\infty, a]$ for some $a$. 
Then the $Tor$-finiteness implies the cut-off the complex
below at a sufficient negative value will be locally free:
the added sheaf will be flat and finitely generated, 
hence locally free.
\epf

We now study the required base change properties.
Let $\psi:\tilde Z\to Z$ be a projective morphism of schemes.
We assume $Z$ has enough locally frees (certainly $Z$
quasi-projective over $\com$ suffices). The functor $\psi^*$ induces
a natural derived functor 
$$L\psi^*:D_{coh}^-(Z)\to D_{coh}^-(\tilde{Z})$$
which sends perfect complexes to perfect complexes.

Let $\phi:\tilde{S}\rarr S $ be a base change of schemes and
consider the Cartesian diagram:
\begin{equation*}
\begin{CD}
\tilde{X}
@>>\tilde{f}>  \tilde{Y}  @>>> \tilde{S}  \\
@V{\phi_X}VV   @V{\phi_Y}VV @V{\phi}VV \\
 X  @>>{f}> Y @>>> S.
\end{CD}
\end{equation*}
In this case, $L\phi_X^*$ and $L\phi_Y^*$ may be defined
on complexes of $S$-flat sheaves by $\phi_X^*$ and
$\phi_Y^*$ respectively.

\begin{lm}
\label{compatt}
For each complex $E\comm \in D_{coh}^-(X)$, there is
a natural isomorphism 
\begin{equation}
\label{bobby}
L\phi_Y^* (Rf_* (E\comm)) \rarr R\tilde{f}_*(L\phi_X^*(E\comm)).
\end{equation}
of complexes in $D_{coh}^-(\tilde{Y})$.
\end{lm}

\bpf
As $f$ is projective, $E\comm$ may be taken to be a complex
of locally free sheaves (bounded from above).
Let ${\mathcal{U}}$ be an $f$-relative Cech covering of
$X$. Then the pull-back covering ${\tilde{\mathcal{U}}}$ is
a $\tilde{f}$-relative Cech covering of $\tilde{X}$ (as may
be checked locally on $Y$). As $E\comm$ is a locally free complex,
$L\phi_X^* (E\comm)$ is just $\phi_X^*(E\comm)$ in $D_{coh}^-(\tilde{X})$.
Hence $R\tilde{f}_*(L\phi_X^* (E\comm))$ is represented by
the simple complex on $\tilde{Y}$ associated to 
\begin{equation}
\label{jfk}
C^p({\tilde{\mathcal{U}}}, \phi_X^* E^q)
\end{equation}
On the other hand, $Rf_*(E\comm)$ is the simple
complex on $Y$ associated to the double complex 
\begin{equation}
\label{rfk}
C^p({\mathcal{U}}, E^q).
\end{equation}
The double complex (\ref{jfk}) is easily seen to be the $\phi_Y$
pull-back of the complex (\ref{rfk}). 
As a consequence, there is a natural map 
\begin{equation}
\label{jonny}
L\phi_Y^*(C^p({\mathcal{U}}, E^q) ) 
\rarr C^p({\tilde{\mathcal{U}}}, \phi_X^* E^q).
\end{equation}
As $X$ is flat over $S$, the complex (\ref{rfk}) is
also $S$-flat. Hence, the map (\ref{jonny}) is
a quasi-isomorphism.
\epf

\subsection{The branch divisor construction}
Let $\omega_{X/S}$ and $\omega_{Y/S}$ denote the
relative dualizing sheaves of the structure maps
$p$ and $q$ respectively.
After constructing a natural perfect torsion complex
$$E\comm=[f^*\omega_{Y/S} \rarr \omega_{X/S}],$$ the branch divisor
is defined by $br(f)= div( Rf_*(E\comm))$ on $Y$.

\begin{lm}
\label{canny}
There is a natural morphism 
$f^*\omega_{Y/S}\rarr \omega_{X/S}$.
\end{lm}
\bpf
The canonical morphism $f^*\Omega_{Y/S}\rarr \Omega_{X/S}$ induces
a morphism 
\begin{equation}
\label{twq}
f^*\omega_{Y/S}=\Lambda^n f^*
\Omega_{Y/S}\to\Lambda^n\Omega_{X/S}.
\end{equation}
Locally on $X$, we have an $S$-embedding $X\rarr M$, 
where $M$ is smooth of relative dimension $n+r$ over $S$ 
and $X$ is a local complete intersection. Let $I=I_{X/M}$. 
There is an exact sequence 
\begin{equation*}
0 \rarr {I/I^2} \rarr {\Omega_{M/S}\otimes \oh_X} \rarr
{\Omega_{X/S}} \rarr 0,
\end{equation*}
where $I/I^2$ and $\Omega_{M/S}\otimes \oh_X$ are locally free sheaves 
on $X$ of ranks $r$ and $n+r$. This sequence yields
 a morphism 
\begin{equation}
\label{zszszs}
\Lambda^n\Omega_{X/S}\otimes 
\Lambda^r(I/I^2)\to \Lambda^{n+r}{\Omega_{M/S}\otimes \oh_X}.
\end{equation}
On the other hand, there is a canonical isomorphism 
\begin{equation}
\label{kkkkk}
\omega_{X/S} \stackrel{\sim}{\rarr} 
Hom(\Lambda^r(I/I^2), \Lambda^{n+r}{\Omega_{M/S}\otimes \oh_X}).
\end{equation}
The morphisms (\ref{zszszs}) and (\ref{kkkkk}) above 
induce a morphism 
\begin{equation}
\label{qwt}
\Lambda^n\Omega_{X/S}\to \omega_{X/S}. 
\end{equation}
It is easily checked the locally defined morphism (\ref{qwt})  
is canonical and hence yields a global morphism on $X$.
The Lemma is established by composing (\ref{twq}) with (\ref{qwt}).
\epf

\begin{lm}
Let
$E\comm =[f^*\omega_{Y/S}\to \omega_{X/S}]$. Then $Rf_*(E\comm)$ 
is a perfect torsion complex in $D^-\coh(Y)$.
\end{lm}
\bpf
Since $E\comm$ is perfect, $Rf_*(E\comm)$ is perfect by Lemma
 \ref{perfperf}. To prove $Rf_*(E\comm)$ is torsion on $Y$,
we may use
Lemmas \ref{qtqtqt} and \ref{compatt} to reduce to the case in which $S$ is
a geometric point.
Then, by property (iii), $X$ is reduced. Let $\nu:\tilde X\rarr X$
be a resolution of singularities.
Let $Z_1$ be the image in $Y$ of the singular locus of $X$, 
and let $Z_2\subset Y$ be the locus where $f\circ \nu$ is not \'etale.
Let $Z=Z_1 \cup Z_2$. As $Rf_*(E\comm)$ is exact on
$Y\setminus Z$, the cohomology of $Rf_*(E\comm)$ is supported 
on $Z$. Since $Z$ is 
a closed subset $Y$ of dimension at most $n-1$, $Z$ does not 
contain any point of 
depth zero (a generic point of a component of $Y$).
\epf

\noindent {\bf Definition.} Let $\br(f)=div(Rf_*(E\comm)))$.
We call $\br(f)$ the
{\em generalized branch divisor} of $f$.

\vspace{+10pt}
\noindent {\bf Base change.}
{\em Let $\phi:\tilde S\to S$ be any morphism. 
Properties (i-iii) hold for $\tilde{X}\rarr 
\tilde{Y} \rarr \tilde{S}$, and}
$$\phi_Y^*(\br(f))=\br(\tilde f).$$

\bpf By Lemma \ref{compatt}, $L\phi_Y^*(Rf_*(E_f\comm))=
R{\tilde{f}}_*
(L\phi_X^*(E_f\comm)) = R{\tilde{f}}_*(E_{\tilde{f}}\comm)$. 
The result then follows from Proposition \ref{mainprop}.v .
\epf

By the base change property, the generalized
branch divisor $\br(f)$ is 
a relative Cartier divisor on $Y$. By {\em relative} we mean
here the restriction to any geometric fiber of $Y$ over $S$
is a well-defined
Cartier divisor.

If $p:X\rarr S$ is smooth and every component of $X$
dominates a component of $Y$, then $\br(f)$  is the
standard branch divisor of $f$.

\section{\bf Stable maps}
\label{exten}
\subsection{Moduli points}
\label{mpoints}
Let $\overline{M}_{g}(D,d)$ be the moduli space of genus $g$,
degree $d>0$ stable maps to a nonsingular curve $D$. 
Let 
\begin{equation}
\label{nnn}
F: {\mathcal{C}} \rarr D \times \overline{M}_{g}(D,d)
\end{equation}
be the universal family of maps over $\overline{M}_{g}(D,d)$.
These objects and morphisms naturally lie in the category of
Deligne-Mumford stacks. We could instead utilize equivariant
constructions in the category of schemes to
study these universal objects (See [FuP], [GrP]). In any case, conditions
(i-iii) of Section \ref{sss} are valid for (\ref{nnn}). Hence, there
exists a relative Cartier divisor $\br(F)$ on
$D \times \overline{M}_{g}(D,d)$ over $D$.

Let $[f:C \rarr D] \in \overline{M}_g(D,d)$ be a moduli
point. We first calculate $\br(f)$ on $D$.
Let $\nu:\tilde{C} \rarr C$ be the normalization map, and
let $\tilde{f}= f \circ \nu$.
Let $N$ be the singular locus of $C$ ($N$ is the union of the
nodal points). 

There are canonical exact sequences
\begin{equation}
\label{redd}
0 \rarr  f^*\omega_D \rarr  \nu_* \tilde{f}^* \omega_D \rarr 
 \oh_N \otimes f^*\omega_D \rarr 0,
\end{equation}
\begin{equation}
\label{bluee}
0 \rarr  \nu_*(\omega_{\tilde C}) \rarr  \omega_C \rarr \oh_N \rarr 0.
\end{equation}
We will use these sequences to express the branch divisor
$br(f)$ as a sum over component contributions.

\begin{lm}
\label{primmm} 
$\br(f)=\br(\tilde f)+2f_{*}(N).$
\end{lm}
\bpf
Since $\nu$ is a finite map,
\begin{equation}
\label{finnn}
R\tilde{f}_{*}([\tilde{f}^*\omega_D
\rarr \omega_{\tilde {C}}]) \stackrel{\sim}{\rarr}
Rf_{*}([\nu_*\tilde f^*\omega_D\rarr \nu_*\omega_{\tilde C}]). 
\end{equation}
Using (\ref{redd}) and (\ref{bluee}), there are a natural distinguished
triangles in $D^-\coh(C)$:
$$[f^*\omega_D \rarr \nu_*\omega_{\tilde C}] \rarr
[\nu_*\tilde f^*\omega_D\rarr \nu_*\omega_{\tilde C}] \rarr
[\oh_N \otimes f^*\omega_D \rarr 0],$$
$$[f^*\omega_D \rarr \nu_*\omega_{\tilde C}] \rarr
[f^*\omega_D\rarr \omega_{\tilde C}] \rarr
[0 \rarr\oh_N].$$
Push-forward
by $Rf_*$ preserves distinguished triangles.
By (\ref{finnn}) and the first triangle,
$$br(\tilde{f})= div Rf_*(
[f^*\omega_D \rarr \nu_*\omega_{\tilde C}]) - f_*(N)$$
(using also properties (iv) and (vi) of Proposition \ref{mainprop}).  
The second triangle yields
$$br(f)= div Rf_*(
[f^*\omega_D \rarr \nu_*\omega_{\tilde C}]) + f_*(N).$$
The Lemma now follows.
\epf

Let $A_1, \ldots, A_a$ be the components of $\tilde{C}$
which dominate $D$, and 
let $B_1, \ldots, B_b$ be
the components of $\tilde{C}$ contracted over $D$.
Let 
$$\{a_i: A_i \rarr D\}, \ \ \{b_j: B_j \rarr D\}$$
denote the natural restrictions of $f$. 
As $a_i$ is a surjective map between nonsingular
curves, the branch divisor $\br(a_i)$ is defined by (\ref{basiccase}). 
Let $b_j(B_j)= p_j \in D$.

\begin{lm} 
\label{seccc}
Let $b: B \rarr p \in D$ be
a contracted component. Then,
$\br(b)=(2g(B)-2)[p]$.
\end{lm}
\bpf
The complex $[b^*\omega_D\rarr \omega_{B}]$ 
is isomorphic to $[\oh_{B} \rarr \omega_{B}]$
with the zero map. By Lemma \ref{smashpoint},
 $div Rf_*(\oh_{B})=\chi(\oh_{B})[p]=(1-g(B))[p]$, 
and $div Rf_*(\omega_{B})=\chi(\omega_{B})[p]=
(g(B)-1)[p]$.
As there is a distinguished triangle in $D^-\coh(B)$:
$$[\oh_{B}] \rarr [\omega_{B}] \rarr
[\oh_{B} \rarr \omega_{B}],$$
we find  $\br(b)= (2g(B)-2) [p]$.
\epf

Lemmas \ref{primmm} and \ref{seccc} prove:
\begin{equation}
\label{ptwise2}
\br(f)= \sum_{i} \br(a_i) +
\sum_j (2g(B_j)-2)[p_j] + 2f_{*}(N).
\end{equation}
The only negative contributions in (\ref{ptwise2}) occur for 
contracted genus 0 components of $\tilde{C}$. However, by stability
such components
must contain at least 3 nodes of $C$.
The Cartier divisor $\br(f)$
is therefore effective for every moduli point $[f:C\rarr D]$.

\subsection{Universal effectivity}
The effectivity of $\br(F)$ over each closed point
of $\overline{M}_{g}(D,d)$ does not guarantee
$\br(F)$ is an effective Cartier divisor on
$D \times \overline{M}_g(D,d)$. The latter
effectivity will now be established.

Let $\pi: {\mathcal{C}} \rarr \overline{M}_{g}(D,d)$ be
the structure map of the universal curve.
There is a canonical exact sequence on $\mathcal{C}$:
\begin{equation}
\label{lsls}
0 \rarr K \rarr F^*\omega_D \rarr \omega_\pi \rarr Q \rarr 0.
\end{equation}
The following vanishing statement will be proven in
Section \ref{vanny}.
\begin{lm}
\label{fitd}
$R^0F_*(K)=0$ and $R^1F_*(Q) =0$.
\end{lm}

Let $E\comm= [F^*\omega_D \rarr \omega_\pi]$ in 
$D^-\coh({\mathcal{C}})$. By definition,
$$br(F)= div(RF_*(E\comm)).$$ The cohomology of $RF_*(E\comm)$
may be computed via a spectral sequence with $E_2$ term:
$$R^1F_*(K) \ \ \ R^1F_*(Q)$$ 
$$R^0F_*(K) \ \ \ R^0F_*(Q)$$
where the grading is -1 for the bottom left corner, 0
for the diagonal, and +1 for the top right corner.
By Lemma \ref{fitd}, we see $RF_*(E\comm)$ has cohomology 
only at grade 0. Hence, locally on $D \times \overline{M}_g(D,d)$,
the complex
$RF_*(E\comm)$ is isomorphic in the derived category to a
finite resolution of the coherent torsion sheaf $H^0(RF_*(E\comm))$.
By Mumford's effectivity result (Proposition 1.ii),
$div(RF_*(E\comm))$ is an effective Cartier divisor on
$D \times \overline{M}_g(D,d)$.

As $br(F)$ is effective {\em and} $\pi$-relatively 
effective, $br(F)$ determines a $\pi$-flat subscheme of
$D\times \overline{M}_{g}(D,d)$. The relative
degree of $\br(f)$ is $r=2g-2-d(2g(D)-2)$.
We have proven:

\begin{tm}
The branch divisor
$\br(F)$ induces a morphism:
$$\gamma: \overline{M}_{g}(D,d) \rarr \text{Hilb}^r(D)= Sym^r(D).$$
\end{tm}

\subsection{Proof of Lemma \ref{fitd}}
\label{vanny}
We follow here the notation of Section \ref{mpoints}.
The first step in the proof is:
\begin{lm}
\label{carey}
The vanishings
$R^0F_*(K)=0$, $R^1F_*(Q)=0$ are equivalent to the
vanishings $R^0\pi_*(K)=0$,
$R^1\pi_*(Q)=0$ respectively.
\end{lm}
\bpf
Let $p: D\times \overline{M}_g(D,d)\rarr \overline{M}_g(D,d)$
be the projection. 
Consider first $K$.
Since $\pi=p\circ F$, there is a spectral
sequences with $E_2$ term:
$$R^1p_*(R^0F_*(K))\ \ \   R^1p_*(R^1F_*(K))$$
$$R^0p_*(R^0F_*(K))\ \ \   R^0p_*(R^1F_*(K))$$
which calculates the sheaves $R^i\pi(K)$.
As both $R^0F_*(K)$ and $R^1F_*(K)$ have support {\em finite}
over $\overline{M}_g(D,d)$, the first row of the above
spectral sequence vanishes. Hence, 
$$R^0\pi_*(K)= R^0p_*(R^0F_*(K)).$$
Moreover as the support of $R^0F_*(K)$ is $p$-finite,
$R^0F_*(K)$ vanishes if and only if $R^0p_*(R^0F_*(K))$
does.
The proof for $Q$ is identical as the supports of the sheaves
$R^iF_*(Q)$ are also $p$-finite.
\epf

\begin{lm}
\label{busey}
$R^0\pi_*(K)=0$.
\end{lm}
\bpf
Let $[f:C \rarr D]$ be a moduli point. Let $[f]\in U\subset
\overline{M}_{g}(D,d)$ where $U$ is an open set.
Let $z\in \Gamma(\pi^{-1}(U),K)$.
The element $z$ is naturally a section of $F^*\omega_D$
over $\pi^{-1}(U)$ which lies in the kernel of
$$F^*\omega_D \rarr \omega_\pi.$$
Let $Spec(A) \subset \overline{M}_g(D,d)$ be any Artinian
local subscheme supported at $[f]$.
We will show the restriction of $z$ to the
closed scheme ${\mathcal{C}}_A= \pi^{-1}(Spec(A))$
vanishes for all such Artinian local subschemes. This vanishing suffices
to prove $z=0$ over a Zariski open neighborhood of $[f]$ by
the Theorem on formal functions (see [Ha1]).

For notational simplicity, let $L= F^*\omega_D$ on ${\mathcal{C}}$.
Let $B\subset C$ be the union of subcurves
contracted by $f$. Since $L|_B$ is trivial, we find the
vanishing condition:
{\em a section of $\Gamma(C, L_{[f]})$ which has
support on $B$ must vanish identically}.

Let $Spec(A) \subset \overline{M}_g(D,d)$ be an Artinian
local subscheme as above. Let $z_A$ be the restriction
of $z$ to ${\mathcal{C}}_A$.
Let $m\subset A$ be the maximal ideal.
We note $z_A$ must have support on $B$ as the sheaf map
$L_A \rarr \omega_{\pi_A}$ is an isomorphism on the
open set $B^c \subset {\mathcal{C}}_A$.
By the flatness of $\pi$, there is an exact sequence
$$0 \rarr m L_A \rarr L_A \rarr L_{[f]} \rarr 0$$
on ${\mathcal{C}}_A$.
By the vanishing condition we see $z_A \in \Gamma({\mathcal{C}}_A,
mL_A)$.
We then use the exact sequences
$$0 \rarr m^{k+1} L_A \rarr m^k L_A \rarr m^k/m^{k+1} \otimes
L_{[f]} \rarr 0$$ to
inductively  prove $z_A \in \Gamma({\mathcal{C}}_A,
m^kL_A)$ for all $k$. Thus $z_A=0$ by the Artinian condition.
\epf

\begin{lm}
\label{ausey}
$R^1\pi_*(Q)=0$.
\end{lm}
\bpf
From sequence (\ref{lsls}), we obtain:
$$R^1\pi_*(F^*\omega_D) \stackrel{i}{\rarr} 
R^1\pi_*(\omega_\pi) \rarr R^1\pi_*(Q) \rarr 0$$
on $\overline{M}_g(D,d)$. It suffices to prove $i$ is a
surjection of sheaves. 
As before, let $[f:C \rarr D]$ be a moduli point. 
Consider the standard
diagram:

\begin{equation*}
\begin{CD}
R^1\pi_*(F^* \omega_D)_{[f]}
@>{i_{[f]}}>>  R^1\pi_*(\omega_\pi)_{[f]}    \\
@V{s}VV   @V{t}VV  \\
 H^1(C, F^*\omega_D)  @>{j}>> H^1(C,\omega_C).
\end{CD}
\end{equation*}
Here, the top line denotes the fiber of the sheaves
at the point $[f]$.
As $R^1 \pi_*(\omega_\pi)$ is locally free
and Serre dual to $R^0\pi_*(\oh_{\mathcal{C}})$ on
$\overline{M}_{g}(D,d)$, the map $t$
is an isomorphism.
As $F^*\omega_D$ is $\pi$-flat, we may apply the cohomology
and base change Theorem (see [Ha1]) to deduce $s$ is surjective.
(As $R^2\pi_*(F^*\omega_D)_{[f]} \rarr H^2(C, F^*\omega_D)$
is trivially surjective and $R^2\pi_*(F^*\omega_D)$
is locally free, the surjectivity of $s$ follows.)
Surjectivity of $i$ locally at $[f]$ is equivalent to
the surjectivity of $i_{[f]}$ by Nakayama's Lemma.
Therefore, the Lemma may be proven by showing $j$ is 
surjective.

It suffices finally to prove $H^1(C,Q)=0$. Again,
let $B\subset C$ be the union of subcurves
contracted by $f$.
Let $I_B \subset \oh_C$ be the ideal sheaf of $B$.
As the map $F^*\omega_D \rarr \omega_C$ is
0 on $B$, we see:
$$Image(F^*\omega_D) \subset I_B \otimes \omega_C.$$
Hence, there is a sequence
$$0 \rarr T \rarr Q \rarr \omega_C|_B \rarr 0$$
where $T$ is easily seen to be a torsion sheaf.
Then, 
$$h^1(C,Q)=h^1(C,\omega_C|_B)=h^0(C, Hom(\omega_C|_B, \omega_C)).$$
The last equality is by Serre duality.
As $B$ is a proper subcurve, the last cohomology group certainly
vanishes.
\epf

Lemmas \ref{carey}-\ref{ausey} combine to prove Lemma \ref{fitd}.
\section{\bf Hurwitz numbers}
\label{hurnum}
\subsection{Integrals}
Let $g\geq 0$ and $d>0$ be integers. Let $b$ be a fixed general
divisor of degree $2g-2+2d$ on $\proj^1$. Let $H_{g,d}$ be the
number of nonsingular genus $g$ curves expressible as $d$
sheeted covers of $\proj^1$ with branch divisor $b$.
There is a long history of the study of $H_{g,d}$
in geometry and combinatorics. The first approach to these
numbers via the combinatorics of the symmetric group was
pursued by Hurwitz in [Hu].

\begin{pr} The Hurwitz numbers are integrals in Gromov-Witten
theory:
\begin{equation*}
H_{g,d}= \int_{ [\overline{M}_{g}(\proj^1,d)]
^{vir}} \gamma^*(\xi^{2g-2+2d}),
\end{equation*}
where $\xi$ is the hyperplane class on 
$\sym^{2g-2+2d}(\proj^1)=\proj^{2g-2+2d}$.
\end{pr}

\bpf We first prove the locus $M_g(\proj^1,d) 
\subset \overline{M}_g(\proj^1,d)$
is nonsingular (of the expected dimension). It suffices to
prove the obstruction space $\text{Obs}(f)$ vanishes.
Let $[f:C \rarr \proj^1]$ be a moduli point with $C$ nonsingular.
There is a canonical right  exact sequence:
$$ H^1(C, T_C) \stackrel{i}{\rarr}
 H^1(C, f^*T_{\proj^1}) \rarr \text{Obs}(f) \rarr 0.$$
Since $d>0$, the sheaf map
$T_C \rarr f^*T_{\proj^1}$ has a torsion quotient.
Hence, $i$ is surjective and $\text{Obs}(f)=0$.
The virtual class of $\overline{M}_{g}(\proj^1,d)$  must then restrict
to the ordinary fundamental class of the open set $M_g(\proj^1,d)$.

Let $r=2g-2+2d$. Let $p_1, \ldots, p_r\in \proj^1$ be distinct points.
By the computation of $\gamma$ on singular curves (\ref{ptwise2}),
we find
$\gamma^{-1}(\sum p_i) \subset M_g(\proj^1,d)$.
By Bertini's Theorem applied to
$$\gamma:M_g(\proj^1,d) \rarr \proj^r,$$
a general divisor $\sum p_i$
intersects the stack $M_{g}(\proj^1,d)$ transversely via $\gamma$ in
a finite number of points. These points are simply the finitely many
Hurwitz covers ramified over $\sum p_i$.
\epf

The first approach to the Hurwitz numbers is via
divisor linear equivalences in $\overline{M}_g(\proj^1,d)$. 
In genus 0 and 1, the divisor of maps $$D_p \subset
\overline{M}_g(\proj^1,d)$$
ramified over a fixed point $p\in \proj^1$ may be expressed
in terms of the boundary divisors of $\overline{M}_g(\proj^1,d)$:
$D_p = \sum_{i} \alpha_i \Delta_i$.
The equation
$$H_{g,d}= D_p \cap \gamma^*(\xi^{r-1})= \sum_i \alpha_i
\Delta_i \cap \gamma^*(\xi^{r-1})$$ then
immediately yields 
recursive relations for $H_{g,d}$:

\begin{equation*}
H_{0,d}= \frac{2d-3}{d} \sum_{i=1}^{d-1} 
\binom{2d-4}{2i-2} i^2(d-i)^2 H_{0,i}H_{0,d-i},
\end{equation*}

\begin{eqnarray*}
H_{1,d} & =&   \frac{d}{6}\binom{d}{2}(2d-1)H_{0,d} \\
&   & +\sum_{i=1}^{d-1} 
\binom{2d-2}{2i-2} (4d-2)i^2(d-i)H_{0,i}H_{1,d-i}.
\end{eqnarray*}

\noindent
The above recursions
were derived by the second author and T. Graber. 
R. Vakil has extended these formulas in genus 0 and 1 by refining the
method
to include non-simple branching. We omit the proofs here
since a uniform treatment may be found in [Va].
Following the
shape of these equations, the recursion 
\begin{eqnarray*}
H_{2,d} & =&  d^2(\frac{97}{136}d-\frac{20}{17})H_{1,d} \\
&   & +\sum_{i=1}^{d-1} 
\binom{2d}{2i-2} (8d-\frac{115}{17}i)i(d-i)H_{0,i}H_{2,d-i}\\
& & + \sum_{i=1}^{d-1} 
\binom{2d}{2i} (\frac{11697}{34}i(d-i) -\frac{3899}{68}d^2)i(d-i)
H_{1,i}H_{1,d-i}
\end{eqnarray*}
was conjectured by the second author and T. Graber in 1997.
Using a completely different combinatorial approach,
Goulden and Jackson have recently proven the genus 2 conjecture
in [GoJ].

The existence of the genus 2 relation does not yet have
a straightforward geometric explanation. In this sense,
it is analogous to the surprisingly simple Virasoro prediction 
for the elliptic Gromov-Witten invariants of $\proj^2$
(see [EHX], [P], [DZ]).
It is not likely such simple recursive formulas occur for
$g\geq 3$.

\subsection{Localization}
Let the torus $\com^*$ act on $ V= \com \oplus \com$
diagonally with weights $[0,1]$ on a basis set $[v_1,v_2]$.
This action induces natural $\gamma$-equivariant
actions on $\overline{M}_g(\proj(V),d)$ and 
$$\proj^r=\text{Sym}(\proj(V))=
\proj(\text{Sym}^rV^*).$$
Moreover, the $\com^*$ action lifts equivariantly to 
the line bundle $$L=\oh_{\proj^r}(1).$$
The choice of equivariant lift to $L$ will be exploited below.
The integral
\begin{equation}
\label{hurww}
H_{g,d}= \int_{ [\overline{M}_{g}(\proj(V),d)]^{vir}} \gamma^*(c_1(L)^r),
\end{equation}
may then be evaluated via 
the virtual localization formula [GrP].

The connected components of the $\com^*$-fixed locus of 
$\overline{M}_{g}(\proj(V),d)$ are indexed by a set of labelled
connected
graphs $\Gamma$ first studied by Kontsevich [Ko].  
The vertices of these graphs lie over the
fixed points $p_1, p_2 \in \proj^1$ and are
labelled with genera (which sum over the graph to $g-h^1(\Gamma)$).
The edges of the graphs lie over $\proj^1$ and
are labelled with degrees (which sum over the 
graph to $d$).
The virtual localization formula of [GrP] yields the equation:
\begin{equation}
\label{ffor}
H_{g,d}= \int_{\overline{M}_g(\proj(V),d)} \gamma^*(c_i(L)^r) =
\sum_{\Gamma} \frac{1}{\text{Aut}(\Gamma)} \int_{\overline{M}_\Gamma} 
\frac{\gamma^*(c_1(L)^r)}{e(N^{vir}_\Gamma)}
\end{equation}
where $\overline{M}_\Gamma$ is a product 
moduli spaces of stable pointed curves and
$\overline{M}_\Gamma/\text{Aut}(\Gamma)$
is the fixed locus associated to $\Gamma$ (see [GrP]).
Moreover, the equivariant Euler class of the virtual
normal bundle, $e(N^{vir}_\Gamma)$,
is explicitly calculated in terms of the
tautological $\psi$ and $\lambda$ classes
on $\overline{M}_\Gamma$. 
Recall the Hodge integrals are the top intersection
products of the $\psi$ and $\lambda$ classes on the
moduli spaces of curves (see [GeP], [FP]).
For each choice of equivariant
lifting of $\com^*$ to $L$, formula (\ref{ffor}) yields
an explicit Hodge integral expression for $H_{g,d}$.

There are exactly $r+1$ distinct
$\com^*$-fixed points of $\proj^r= \proj(\text{Sym}^rV^*)$.
For $0\leq a \leq r$,
Let $p_{a}$ denote the fixed point $v_1^{*a} v_2^{*(r-a)}$. 
The canonical $\com^*$-linearization on $L= \oh(1)$
has weight $w_{a}= r-a$ at $p_{a}$.
Let $L_i$ denote the unique 
$\com^*$-linearization of $L$
satisfying $w_{i}=0.$
We note the weight at $p_0$ of $L_i$ is $i$.
We may rewrite (\ref{ffor}) as:
\begin{equation}
\label{fforr}
H_{g,d}= 
\sum_{\Gamma} \frac{1}{\text{Aut}(\Gamma)} \int_{\overline{M}_\Gamma} 
\frac{\prod_{i=1}^{r} \gamma^*(c_1(L_i))}
{e(N^{vir}_\Gamma)}.
\end{equation}
This choice of linearization for the integrand will
lead to the simplest localization formula.

The morphism $\gamma$ associates a unique fixed point $p_\Gamma$ to
each graph $\Gamma$:
$$ \gamma(\overline{M}_\Gamma/\text{Aut}(\Gamma))= p_\Gamma.$$
Let $p_\Gamma=p_i\neq p_0$.
Then, $\gamma(L_i)$ is a trivial bundle with trivial linearization
when restricted to the fixed locus $\overline{M}_\Gamma/\text{Aut}(\Gamma)$.
The $\Gamma$-contribution to
the sum (\ref{fforr}) therefore vanishes.
We must only consider those graphs $\Gamma$ satisfying
$p_\Gamma=p_0$ in the sum (\ref{fforr}).

The point $p_0= [v_2^{*r}]$ corresponds to the
divisor $r[v_1]$ on $\proj(V)$. It is a very
strong condition for a stable map $[f:C \rarr \proj(V)]$ to have
$br(f)$ supported at the single point $[v_1]$ -- all nodes, collapsed
components, and ramifications must lie over $[v_1]$.
Hence, if $p_\Gamma=p_0$, the graph $\Gamma$ may not have
any vertices of positive genus or valence greater than 1 lying
over $[v_2]$. Moreover, the degrees of the edges of $\Gamma$
must all be 1. Exactly one
graph $\Gamma_0$ satisfies these conditions:
$\Gamma_0$ has a single genus $g$ vertex lying over $[v_1]$ which
is incident to exactly
$d$ degree 1 edges (the vertices over $[v_2]$ are
degenerate of genus 0).

The sum (\ref{fforr}) contains only one term:
\begin{equation}
\label{fforrr} 
H_{g,d}= 
\frac{1}{\text{Aut}(\Gamma_0)} \int_{\overline{M}_{\Gamma_0}} 
\frac{\prod_{i=1}^{r} \gamma^*(c_1(L_i))}
{e(N^{vir}_{\Gamma_0})}.
\end{equation}
By definition (see [GrP]),
$\overline{M}_{\Gamma_0}= \overline{M}_{g,d}$.
Since the automorphism group of $\Gamma_0$
is the full permutation group of the edges, $\text{Aut}(\Gamma_0)= d!$.
The virtual normal bundle contribution is calculated in
[GrP]:
$$\frac{1}{e(N^{vir}_{\Gamma_0})}= 
\frac{1-\lambda_1+\lambda_2 -\lambda_3 + \ldots + (-1)^g \lambda_g}
{\prod_{i=1}^d (1-\psi_i)}.$$
Finally, the integrand $\prod_{i=1}^r \gamma^*(c_1(L_i))$
restricts to a term of pure weight $r!$.
We have proven:

\begin{tm}
\label{hodgehurr}
$$H_{g,d}= \frac{(2g-2+2d)!}{d!} \int_{\overline{M}_{g,d}}
\frac{1-\lambda_1+\lambda_2 -\lambda_3 + \ldots + (-1)^g \lambda_g}
{\prod_{i=1}^d (1-\psi_i)},$$
for $(g,d) \neq (0,1), (0,2)$.
\end{tm}

In genus 0, there is a  well-known formula for
the $\psi$ integrals:
$$\int_{\overline{M}_{0,n}} \psi_1^{a_1} \cdots \psi_n^{a_n} =
\binom{n-3}{a_1,\ldots, a_n}$$
(see [W]).
The genus 0 formula
\begin{equation}
\label{gg00}
H_{0,d}= \frac{(2d-2)!}{d!} d^{d-3}
\end{equation}
then follows immediately from Theorem 2.
Equation (\ref{gg00}) was first found by Hurwitz.

\vspace{+10pt}
\noindent Department of Mathematics \\
\noindent University of Trento \\
\noindent fantechi@alpha.science.unitn.it

\vspace{+10 pt}
\noindent
Department of Mathematics \\
\noindent California Institute of Technology \\
\noindent rahulp@cco.caltech.edu

\end{document}